\documentstyle{article}
\begin{document}

\newtheorem{satz}{Theorem}[section]
\newtheorem{lem}[satz]{Lemma}
\newtheorem{prop}[satz]{Proposition}
\newtheorem{kor}[satz]{Corollary}
\noindent
RINGS OF GLOBAL SECTIONS IN TWO\--DIMENSIONAL SCHEMES
\par
\bigskip
\noindent
\centerline{Holger Brenner, Bochum}
\par
\bigskip
\noindent
\begin{abstract}
In this paper we study the ring of global sections $\Gamma(U,{\cal O})$
of an open subset
$U=D(I) \subseteq {\rm Spec}\,A$, where $A$ is a two-dimensional
noetherian ring. The main concern is to give a geometric criterion when
these rings are
finitely generated, in order to correct an invalid statement of Schenzel
in \cite{schenzel88}.
\end{abstract}
\par\bigskip
\noindent
\section{Introduction}
\par
\bigskip
\noindent
Let $A$ be a noetherian ring with an ideal $I \subseteq A$ and
$U=D(I) \subseteq {\rm Spec}\, A$ the corresponding open subset.
If $U$ is an affine scheme, then the ring of global sections
$B=\Gamma(U,{\cal O}_X)$ --which is also called the ideal-transform
$T(I)$-- is of finite type over $A$.
The converse is by no means true, in dimension two however
we have the following result due to Eakin et. al.
(\cite{eakin}, Theorem 3.2):
Suppose $A$ is a local excellent
\footnote{In fact the result was stated under the somewhat weaker conditions
that the normalization is finite and the local rings of the normalization
are analytically irreducible, instead of excellent.}
Cohen-Macaulay domain of dimension
two, and let $I$ be an ideal of height one.
Then (among other characterizations)
$D(I)$ is affine if and only if $B$ is noetherian if and only if
$B$ is of finite type over $A$.
\par\smallskip\noindent
Schenzel states in \cite{schenzel88}, Theorem 4.1 and 4.2, that
this holds more general for two-dimensional excellent local domains.
However, this is not true, as the following example shows.
\footnote{The mistake in \cite{schenzel88} is at the end of the proof of
theorem 4.1, where the statement $T \subseteq T_N$ is wrong.}
\par
\bigskip
\noindent
{\bf Example.}
Let $X={\rm Spec}\, A$ be an affine excellent irreducible surface
which is regular outside
one single closed point $P$ and such that in the normalization
two points $P_1,P_2$ lie over $P$. Outside these points the normalization
mapping $\tilde{X} \longrightarrow X$ is isomorphic.
\par\smallskip\noindent
Let $Y=V(I)$ be the image of an irreducible curve $Y'$ passing through $P_1$,
but not through $P_2$.
Then $U=X-Y$ is not affine, since the preimage of $Y$ consists of the curve
$Y'$ and of the isolated point $P_2$.
On the other hand, $U=X-Y$ is normal and isomorphic to
$\tilde{X}-Y'-P_2 $, so the rings of global sections are identical. Since
$\tilde{X}$ is normal, this ring equals also the ring of global sections of
$\tilde{X}-Y'$. $\tilde{X}$ is a normal excellent affine surface, thus
the complement of a curve is affine, and $B$ is finitely generated.
For an explicit example see below.
\par
\bigskip
\noindent
In this paper we give a criterion for two-dimensional local rings to
decide the finiteness of the ring of global sections of $U=D(I)$, $I$
an ideal of height one.
The criterion is based on the combinatoric of the components in
the completion $\hat{A}$ of $A$. It says that in case $U$ is not affine
the ring of global sections of 
$U$ is not finitely generated if and only if there exists an irreducible
component of ${\rm Spec}\, \hat{A}$
where $U$ is affine and a component where $U$ is not affine
such that their intersection is one-dimensional.
\par\smallskip\noindent
The criterion is (due to the connectedness theorem of Hartshorne)
seen to be fulfilled in case $A$ is Cohen-Macaulay,
thus we recover the result of Eakin et. al. as a corollary (Cor. 2.4).
Another consequence is that if $D(I)$ is non-affine and connected, then
$T(I)$ is not noetherian (Cor. 2.3).
\par\smallskip\noindent
In the third section we extend the result to the non-complete case
and describe the conditions used in the criterion in terms
of the normalization.
\par
\bigskip
\noindent
\section{The complete case}
\par
\bigskip
\noindent
Let $X={\rm Spec}\, A$ be the spectrum
of a local complete noetherian ring $A$ of dimension $2$, and
let $P$ denote the closed point.
Let $X_j=V({\bf p}_j)={\rm Spec}\, A/{{\bf p}_j}$ be the irreducible components of $X$
corresponding to the minimal primes ${\bf p}_j,\, j \in J$.
\par\smallskip\noindent
Let $I$ be an ideal in $A$, $Y=V(I)$
and $U=D(I)$.
$U$ is affine if and only if $U_j=U \cap X_j$ is affine on every
component, and this is due to the theorem of Lichtenbaum-Hartshorne
(see \cite{brodman}, 8.2.1) the case if and only if
${\rm ht}\, I(A/{{\bf p}_j}) \leq 1$ for every $j \in J$.
Thus $U$ is not affine if and only if there exists a two-dimensional
component $X_j$ where $Y_j=Y \cap X_j$ consists just of the single point $P$.
\par\smallskip\noindent
We want to know for an ideal $I$ of height one whether the ring
of global sections of $D(I)$ is finitely generated.
If $D(I)$ is affine, this is the case, so we suppose
furtheron that $D(I)$ is not affine.
We divide $J=J_0 \cup J_1$ in such a way, that for
$ j \in J_1$ the open subsets $U_j  \subseteq X_j$
are affine and for $ j \in J_0$ not.
Thus the $X_j, \, j \in J_0$, are the two-dimensional components of $X$
where $Y_j$ is just the closed point. The affineness of $U$ is equivalent
with $J_0 = \emptyset$.
\par\smallskip
\noindent
Put ${\bf a}_0= \bigcap_{j \in J_0} {\bf p}_j$ and
${\bf a}_1= \bigcap_{j \in J_1} {\bf p}_j$ and
$X_0={\rm Spec}\, A/ {\bf a}_0$, $X_1={\rm Spec}\, A/{\bf a}_1$.
We denote the structure sheaves on these closed subschemes
of $X$ with ${\cal O}_i,\, i=0,1$.
\par\smallskip\noindent
Furthermore we put $U_i=U \cap X_i, i=0,1$, considered as an open subset
in $X_i$ with the induced scheme structure, put $B_i=\Gamma(U_i,{\cal O}_i)$.
$U_1={\rm Spec}\, B_1 \subset X_1$ is affine, $U_0$ is not affine.
The closed embedding $X_1 \hookrightarrow X$ yields a
(closed) restriction map
$\Gamma(U,{\cal O}_X) \longrightarrow \Gamma(U_1,{\cal O}_1)=B_1$.
\par\smallskip\noindent
Finally, let ${\bf b}={\bf a}_0+{\bf a}_1 \subseteq {\bf m}$ and $R=A/{\bf b}$.
$R$ is a zero- or one-dimensional local complete noetherian ring, let
$Z={\rm Spec}\, R$ und $Z^\times =D({\bf m}) \subset Z$.
The dimension of $Z=V({\bf b})=X_0 \cap X_1$ is the crucial point for
$\Gamma(U,{\cal O}_X)$ to be noetherian or not.
\par\noindent
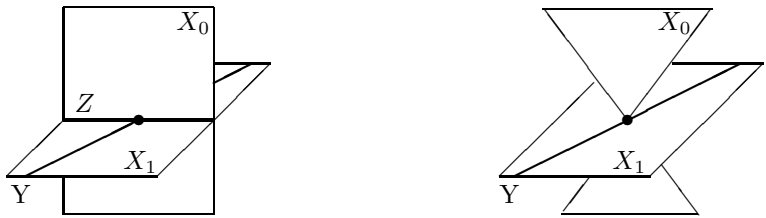
\unitlength0.5cm
\begin{picture}(30,8.5)
\put(5,4.5){\circle*{0.25}}
\put(3.3,4.7){$Z$}
\put(1.5,3){\line(1,1){1.5}}
\put(5.5,3){\line(1,1){3}}
\thicklines
\put(3,4.5){\line(1,0){4}}
\put(2,3){\line(2,1){3}}
\put(7,5.5){\line(2,1){1}}
\put(1.6,2.3){Y}
\thinlines
\put(3,4.5){\line(0,1){3}}
\put(3,7.5){\line(1,0){4}}
\put(7,2){\line(0,1){5.5}}
\put(3,2){\line(1,0){4}}
\put(3,2){\line(0,1){1}}
\put(1.5,3){\line(1,0){4}}
\put(7,6){\line(1,0){1.5}}
\put(4.6,3.2){$X_1$}
\put(6, 6.9){$X_0$}

\put(18, 4.5){\circle*{0.25}}
\put(14.6, 3){\line(1,1){2.5}}
\put(18.6,3){\line(1,1){3}}
\put(19.2,6){\line(1,0){2.4}}
\put(14.6,3){\line(1,0){4}}
\thicklines
\put(15,3){\line(2,1){6}}
\put(14.6,2.3){Y}
\thinlines

\put(17.6,3.2){$X_1$}
\put(18.8,6.9){$X_0$}

\put(18,4.5){\line(3,4){2.2}}
\put(18,4.5){\line(-3,4){2.2}}
\put(15.75, 7.45){\line(1,0){4.5}}

\put(18.9, 3.3){\line(3,-4){1.0}}
\put(17, 3){\line(-3,-4){0.7}}
\put(16.25, 2){\line(1,0){3.6}}

\end{picture}
\par
\noindent
For our proof we have to put on $A$ the condition $S_1$ of Serre,
meaning that every associated prime of $A$ is minimal,
equivalently that every zero-divisor lies in a 
minimal prime or that every ideal of height one contains a non-zero-divisor.
This is fulfilled for example if $A$ is reduced.
\begin{satz}\hspace{-0.5em}{\bf .}
Let $A$ be a two-dimensional complete local noetherian ring,
fulfilling the condition $S_1$.
Let $I$ be an ideal of height one and suppose that $U=D(I)$ is not affine.
Then the following are equivalent.
\par\smallskip\noindent
{\rm (1)} $\Gamma(U,{\cal O}_X)$ is not of finite type.
\par\smallskip\noindent
{\rm (2)} $\Gamma(U,{\cal O}_X)$ is not noetherian.
\par\smallskip\noindent
{\rm (3)} The image of
$\Gamma(U,{\cal O}_X) \longrightarrow \Gamma (U_1,{\cal O}_1)$
is not noetherian.
\par\smallskip\noindent
{\rm (4)} The intersection $Z$ of $X_0$ and $X_1$ is one-dimensional.
\end{satz}
{\it Proof}.
The implications $(3) \Rightarrow (2)$ and $(2) \Rightarrow (1)$ are clear.
$(1) \Rightarrow (4)$. Suppose $Z =\{P\} $ is only the closed point.
Then $U$ is the disjoint union
of $U_0$ and $U_1$ (both closed hence open in $U$). Thus we have
$$\Gamma(U,{\cal O}_X)
= \Gamma(U_0,{\cal O}_0) \oplus \Gamma(U_1,{\cal O}_1)\, .$$
Since $U_1$ is affine, the second component is of finite type.
Since $U_0=X_0 -\{P\}$, the mapping
$A/{\bf a}_0 \longrightarrow \Gamma(U_0,{\cal O}_X)$ is also of
finite type, see Lemma 2.2 (1).
\par\smallskip\noindent
So we have to show $(4) \Rightarrow (3)$. We denote the image of
$\Gamma(U,{\cal O}_X) \longrightarrow \Gamma(U_1,{\cal O}_1)$ by $C$.
\par\smallskip\noindent
Let $h \in A$ be an element such that in $Z={\rm Spec}\, R$
we have $V(h)=V({\bf m})=\{P\}$.
Thus $1/h$ is a function defined on $Z^\times = Z-\{P\} =D(h)$.
Since $Z^\times \hookrightarrow U_1$ is a closed embedding
and since $U_1$ is affine, there exists a function
$q \in B_1=\Gamma(U_1,{\cal O}_1)$ with $q\mid_{Z^\times}=1/h$.
\par
\noindent
\unitlength0.4cm
\begin{picture}(20,9)

\put(3,4){$U_0$}
\put(3,7){$X_0$}
\put(7,4){$Z^\times$}
\put(7,7){$Z$}
\put(11,4){$U_1$}
\put(11,7){$X_1$}

\put(9,4){$\hookrightarrow$}
\put(9,7.){$\hookrightarrow$}
\put(5,4.){$\hookleftarrow$}
\put(5,7.){$\hookleftarrow$}

\put(3.1,5.5){$\cup $}
\put(7.1,5.5){$\cup $}
\put(11.1,5.5){$\cup $}

\put(7,1){$U$}
\put(7.1,2.5){$\cap$}
\put(3.8,3.5){\vector(1,-1){2}}
\put(10.5,3.5){\vector(-1,-1){2}}
\end{picture}
\par\noindent
Let $a \in {\bf b} \subset A$ be a regular element (i.e. a non-zero-divisor)
inside the describing ideal of $Z$.
The functions $aq^n$ are defined on $U_1$ and the restrictions to
$Z^\times$ are zero,
thus they are extendible to $Z$.
Since $Z \hookrightarrow X_0$ is closed and $X_0$ is affine, these
functions are also extendible to $X_0$ and in particular to $U_0$.
So we may assume that these functions are defined on $U$ and we
see that they lie in $C$.
\par\smallskip\noindent
Consider in $C$ the ideal $(a,aq,aq^2,aq^3,...)$ spanned by this functions,
and suppose that it is finitely generated.
Then we have an equation
$$aq^{n+1}= a_n a q^n +...+ a_1 aq +a_0a $$
with $a_i \in C \subset B_1$.
We may assume that $a_i \in \Gamma(U,{\cal O}_X)$.
Since $a$ is regular in $A$, it is
also a regular in $A/{\bf a}_i$.
(For if $ax \in {\bf a}_i=\bigcap_{j \in J_i}{\bf p}_j$,
we have $ax \in {\bf p}_j$ for all $j \in J_i$ and thus
$x \in {\bf p}_j$ for all $j \in J_i$, so $x=0$ mod ${\bf a}_i$.)
Since the restriction
$A/{\bf a}_1=\Gamma(X_1,{\cal O}_1) \longrightarrow \Gamma(U_1,{\cal O}_1)$
is injective, $a$ is also a regular element in $B_1$.
\par\smallskip\noindent
This yields in $B_1$ (on $U_1$) the equation
$q^{n+1}= a_n q^n+...+a_1q +a_0$.
This equation restricted to $Z^\times \subseteq U_1$
yields an integral equation for $q=1/h$ over $R[a_i'] \subseteq R_h$,
where the $a_i'$ denote the restrictions of $a_i$ on
$R_h=\Gamma(Z^\times,{\cal O}_Z)$.
\par\smallskip\noindent
We claim that the $a_i'$ are integral over $R$:
Consider the elements $a_i \in \Gamma(U,{\cal O}_X)$ as functions
on $U_0$ --as elements of $B_0$.
Since $U_0=X_0 - \{P\}$, the $a_i \in B_0$ are integral over
$A/{\bf a}_0=\Gamma(X_0,{\cal O}_0)$, see Lemma 2.2.
The closed embeddings
$(Z^\times \subset Z) \hookrightarrow (U_0 \subset X_0)$
show that the $a_i'$ are integral over $R=\Gamma(Z,{\cal O}_Z)$.
It follows that $q\mid_Z=1/h$ would be integral over $R$, but this is
not possible. \hfill $\Box$
\par
\noindent
\begin{lem}\hspace{-0.5em}{\bf .}
Let $A$ be a local noetherian ring of dimension two fulfilling $S_1$.
Let $\bf m$ be the maximal ideal and $B=\Gamma(D({\bf m}),{\cal O})$ the ring
of global sections. Then the following hold.
\par\smallskip\noindent
{\rm (1)} $A \longrightarrow B$ is of finite type.
\par\smallskip\noindent
{\rm (2)} If furthermore all components of ${\rm Spec}\, A$ have dimension two,
$B$ is even finite over $A$.
\end{lem}
{\it Proof}.
We first proof the second part, using \cite{EGAIV}, 5.11.4
(or \cite{binsto1}, 2.5.).
A point $x \in {\rm Ass}\, {\cal O}_X$ has height zero, for every ideal
of bigger height contains a regular element. The closure $\bar{x}$ is
a two-dimensional component and therefore the point $P$ has codimension
two on it.
\par\smallskip\noindent
The first part follows from the second part.
The one-dimensional components of $X$ meet the other components
only in the closed point, thus the punctured curves are
connected components of $W=D({\bf m})$.
These are affine and of finite type. \hfill $\Box$
\par
\bigskip
\noindent
We deduce from the theorem two corollaries.
\begin{kor}\hspace{-0.5em}{\bf .}
Let $A$ be a local complete noetherian ring of dimension two
fulfilling $S_1$.
Let $I$ be an ideal of height one.
If $U=D(I)$ is connected and $\Gamma(U,{\cal O}_X)$ is of finite type,
then $U$ is affine.
\end{kor}
{\it Proof}.
Suppose $U$ is not affine, then in the partition described above
$U_0$ is not empty, and $U_1$ is not empty since $I$ is of height one.
Put $Z=X_0 \cap X_1$. Since $U$ is connected,
$U_0$ and $U_1$ are not disjoint, thus $Z$ does not consist only of the
closed point, it must be a curve. Then due to the theorem
the ring of global sections
can not be noetherian. \hfill $\Box$
\par
\bigskip
\noindent
We recover the result of Eakin et. al. in the Cohen-Macaulay case.
\begin{kor}\hspace{-0.5em}{\bf .}
Let $A$ be a local complete noetherian Cohen-Macaulay ring of dimension two.
Let $I$ be an ideal of height one. Then $U=D(I)$ is affine if and only
if its ring of global sections is of finite type {\rm (}or noetherian{\rm )}.
\end{kor}
{\it Proof}.
Again, suppose $U$ to be not affine, put $X=X_0 \cup X_1$ as before with the
describing ideals ${\bf a}_0$ and ${\bf a}_1$.
Then ${\bf a}_0 \cap {\bf a}_1$ is nilpotent,
thus due to the connectedness theorem of Hartshorne
(see \cite{eisenbud}, theorem 18.12) the ideal
${\bf a}_0 +{\bf a}_1$ has height one. Since it describes the intersection,
$Z=X_0 \cap X_1$ is one-dimensional and $\Gamma(U,{\cal O}_X)$
is not noetherian.  \hfill $\Box$
\par
\bigskip
\noindent
{\bf Example.}
Of course, $U=D(I)$ can be affine without being connected.
$A=K[[x,y,z]]/(xy)$ is Cohen-Macaulay ($K$ a field),
the complement of the common
axis $V(x,y)$ is affine, but not connected.
\par
\bigskip
\noindent
{\bf Remark.}
We may associate to a complete local ring of dimension two a graph $\Gamma$
in such a way, that for each irreducible two-dimensional component
we associate a
point, and two points are connected by an edge if and only if the intersection
of the corresponding components is one-dimensional.
Then an open subset as above yields a partition
$\Gamma =\Gamma_0 \cup \Gamma_1$, and the ring of global sections is
noetherian if and only if there is no edge between points of $\Gamma_0$
and of $\Gamma_1$.
\par
\bigskip
\noindent
\section{Interpretation in the Normalization}
\par
\bigskip
\noindent
We want to extend the result from the complete case to the general
case.
Suppose we are given a curve $V(I) \subseteq {\rm Spec}\, A$
where $A$ is a two-dimensional noetherian domain. Then $\Gamma(D(I),{\cal O})$
is of finite type if this is true in every (closed) point
$x \in {\rm Spec}\, A$, see \cite{bingener}.
Furthermore, we have the following Lemma.
\begin{lem}\hspace{-0.5em}{\bf .}
Let $A \longrightarrow A'$ be faithfully flat and let
$U \subseteq {\rm Spec}\, A$ denote an open subset
with preimage $U'$. Then $B=\Gamma(U,{\cal O})$ is of finite type over $A$
if and only
if $B'=\Gamma(U',{\cal O}')$ is of finite type over $A'$.
\end{lem}
{\it Proof}. 
We have $B' = B \otimes_A A'$ due to flatness.
This yields the first implication. If $B'$ is of finite type,
we may assume that it is generated by finitely many elements of $B$, thus 
there is a surjection $A'[T_1,...,T_n] \longrightarrow B'=B\otimes_AA'$
induced by $A[T_1,...,T_n] \longrightarrow B$. Due to faithfull, this must also
be surjective. \hfill $\Box$
\par
\bigskip
\noindent
Therefore the condition in the theorem
that $\Gamma(U,{\cal O})$ is of finite type is
preserved by passing to the completion, and we may skip in Cor. $2.4$
the assumption of completeness.
\par\smallskip\noindent
So we take a look at the condition that the intersection of two components in
the completion is one-dimensional, and we want to describe it in terms of
the normalization of $A$.
For this we recall some correspondences between normalization and completion,
see \cite{EGAIV}, 7.6.1 and 7.6.2.
Let $X$ be the spectrum of a local excellent domain $A$ with completion
$\hat{X}$ and normalization $\tilde{X}$.
Then the normalization of $\hat{X}$ equals the completion of
$\tilde{X}$ (semilocal), and this consists of connected components being
the normalizations of the irreducible components of $\hat{X}$
and the completion of the localizations of $\tilde{X}$ as well.
In particular, there is a correspondence between the irreducible
components of $\hat{X}$ and the closed points of $\tilde{X}$.
\par\smallskip\noindent
For a closed subset $C \subseteq X$ the completion of $C$ equals the preimage
of $C$ in $\hat{X}$ yielding a canonical inclusion $\hat{C} \subseteq \hat{X}$.
The irreducible components of $\hat{C}$ correspond again to closed points
of $\tilde{C}$, but this is of course not the preimage of $C$ in the
normalization $\tilde{X}$.
\begin{lem}\hspace{-0.5em}{\bf .}
Let $A$ be an excellent local domain of dimension two,
$P_0 \in \tilde{X}$ the closed point
on $\tilde{X}$ corresponding to
the irreducible component $X_0$ of the completion $\hat{X}$.
Let $C \subset X$ be an irreducible curve and let
$D \subset \tilde{X}$ be the preimage of $C$ without the isolated points.
\par\smallskip\noindent
{\rm (1)}
There exists
an irreducible component of $\hat{C}$ on $X_0$ if and only if
$P_0$ is not an isolated point on $\varphi^{-1}(C)$
{\rm (}$\varphi: \tilde{X} \longrightarrow X$ normalization map {\rm )}.
\par\smallskip\noindent
{\rm (2)}
The irreducible component $C_0$ of $\hat{C}$ lies on $X_0$
if and only if there exists a point $R \in \tilde{D}$ over $P_0$
mapping to the point $Q_0 \in \tilde{C}$ corresponding to $C_0$.
\par\smallskip\noindent
{\rm (3)}
The component $C_0$ of $\hat{C}$ connects the irreducible components $X_1$
and $X_2$ of $\hat{X}$ if and only if the corresponding point
$Q_0 \in \tilde{C}$
is reached by points $R_1,R_2 \in \tilde{D}$ lying over $P_1$ and $P_2$.
\end{lem}
{\it Proof.}
(1).
We consider the mapping (completion)
$\tilde{X_0} \longrightarrow \tilde{X}_{P_0}$, where
$\tilde{X}_{P_0}$ means the localization at $P_0$.
The preimage of $C \subset X$ in $\tilde{X}_{P_0}$ is just the closed
point if and only if this is true in $\tilde{X_0}$, and this is the case
if and only if $\hat{C}$ is zero-dimensional on $X_0$.
\par\smallskip\noindent
(2).
The preimage of $\hat{C}$ in $\tilde{\hat{X}}$ without the isolated points
equals $\hat{D}$, being the preimage of $D$.
The statement $C_0 \subset X_0$ is equivalent to the statement that
there exists an irreducible component
$D_0 \subseteq \hat{D} \subset \tilde{\hat{X}}$
over $C_0$ lying on $\tilde{X_0}$.
Let $R$ be the point on $\tilde{D}$ corresponding to the component
$D_0 \subseteq \hat{D}$.
Suitable diagramms show that
$D_0$ dominates $C_0$ is equivalent with $R$ maps to $Q_0$ and that
$D_0 \subseteq \tilde{X_0}$ is equivalent with
$R$ maps to $P_0$
\par\smallskip\noindent
(3) follows from (2).         \hfill $\Box$
\par
\bigskip
\noindent
This motivates the following definition.
\par
\bigskip
\noindent
{\bf Definition.}
Let $X$ denote a reduced irreducible noetherian scheme,
$\varphi: \tilde{X} \longrightarrow X$ its normalization,
$P \in X$ a closed point and $P_1,P_2 \in \tilde{X}$,
$\varphi(P_1)=\varphi(P_2)=P$.
We call an irreducible curve $C \subset X$ a
\emph{melting curve}
for the points
$P_1$ and $P_2$ if and only if $P_1,P_2$ are not isolated
on $\varphi^{-1}(C)$ and there exist points
$R_1,R_2 \in \tilde{D}$ ($D$ as in lemma 3.2)
over $P_1,P_2$ mapping to one
common point $Q \in \tilde{C}$.
\begin{satz}\hspace{-0.5em}{\bf .}
Let $X={\rm Spec}\, A$, where $A$ is an excellent local domain of dimension
two.
Then the intersection of the components $X_1$ and $X_2$ on $\hat{X}$
is one-dimensional if and only if there exists a melting curve for
$P_1,P_2 \in \tilde{X}$.
\end{satz}
{\it Proof}.
If $C$ is a melting curve for $P_1$ and $P_2$ with common point
$Q$ as in the definition, then the previous proposition
says that the corresponding component
$C_0$ lies on $X_1$ and $X_2$, thus the intersection is one-dimensional.
\par\smallskip\noindent
For the converse, let $C_0$ be an irreducible curve on $X_1 \cap X_2$
with prime ideal ${\bf q} \subset \hat{A} $ of height one.
Then ${\bf p}={\bf q} \cap A$ is also of height one.
For $\bf q$ is not a normal point of $\hat{A}$, since on the normalization
there are at least two points above it. 
Then also $\bf p$ is not a normal point, because
the normal locus commutes with completion under the condition of
excellence (see \cite{EGAIV}, 7.8.3.1.)
Thus ${\rm ht}\, {\bf p} =1$, $C=V({\bf p})$ is a curve,
$C_0$ a component of its completion
and we may apply the previous proposition. \hfill $\Box$
\par
\noindent
\begin{prop}\hspace{-0.5em}{\bf .}
Let $P_1,P_2$ be two closed points in the normalization $\tilde{X}$
over $P \in {\rm Spec}\, A$, where $A$ is a two-dimensional noetherian
domain. Then the following hold.
\par\smallskip\noindent
{\rm (1)}
If there exist two different irreducible curves $C_1,C_2$
with $P_i \in C_i=V({\bf q}_i)$ on $\tilde{X}$
such that ${\bf q}_1 \cap A={\bf q}_2 \cap A={\bf p}$, then $C=V({\bf p})$
is a melting curve for $P_1,P_2$.
\par\smallskip\noindent
{\rm (2)}
If $C$ is normal {\rm (}or analytically irreducible{\rm )}
and $P_1$ and $P_2$ are not isolated on $\varphi^{-1}(C)$,
then $C$ is a melting curve.
\par\smallskip\noindent
{\rm (3)}
If $P_1,P_2 \in C'$ is irreducible and $\varphi(C')=C$ is a melting
curve, then $\varphi|_{C'}: C' \longrightarrow C$ is not birational.
A melting curve lies in the non-normal locus.
\end{prop}
{\it Proof}.
(1) Both mappings $C_1 \longrightarrow C$ and $C_2 \longrightarrow C$
are surjective, and this is then also true for the normalizations.
Thus for any closed point $Q \in \tilde{C}$ there are points on
$\tilde{C_i}$ over $P_i$ mapping to $Q$.
\par\smallskip\noindent
(2) If $C$ is analytically irreducible, then any closed point of $\tilde{D}$
maps to the only closed point of $\tilde{C}$.
\par\smallskip\noindent
(3) Suppose $C' \longrightarrow C $ is birational.
Then the normalizations of these curves are the same, and
different points cannot be identified.
If the generic point of a curve $C$ is normal, then $D$ consists just
of one irreducible component, and $D \longrightarrow C$ is
birational.  \hfill $\Box$
\par
\bigskip
\noindent
{\bf Examples.}
We give some typical examples of (non-)melting curves to illustrate the
cases the previous proposition is talking about.
They are given by mappings ${\bf A}_K^2 \longrightarrow {\bf A}_K^n$ such that
the affine plane is the normalization of the image ($K$ is a field).
\par\smallskip\noindent
(1) $(x,y) \longmapsto (x,y^3-y,y^2-1)$.
This identifies the two different curves $V(y-1)$ and $V(y+1)$.
The common image curve $C$ is a melting curve.
\par
\noindent
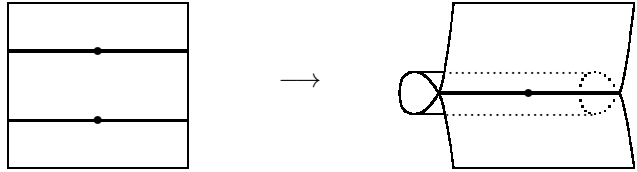
\unitlength0.4cm
\begin{picture}(30,8.5)
\put(5, 3.6){\circle*{0.25}}
\put(5, 5.9){\circle*{0.25}}

\put(2,2){\line(1,0){6}}
\put(2,2){\line(0,1){5.5}}
\put(2,7.5){\line(1,0){6}}
\put(8,2){\line(0,1){5.5}}

\thicklines
\put(2,3.6){\line(1,0){6}}
\put(2,5.9){\line(1,0){6}}
\thinlines

\put(11, 4.7){$\longrightarrow$}

\put(19.3, 4.5){\circle*{0.25}}

\put(16.8, 7.5){\line(1,0){6}}
\put(16.8, 2){\line(1,0){6}}
\thicklines
\put(16.3, 4.5){\line(1,0){6}}
\thinlines
\put(15.5,3.8){\line(1,0){1.04}}
\put(15.5,5.2){\line(1,0){1.0}}
\put(16.5, 3.89){$_{\ldots \ldots \ldots \ldots \ldots \ldots \ldots \ldots}$}
\put(16.4, 5.27){$_{\ldots \ldots \ldots \ldots \ldots \ldots \ldots \ldots}$}
\linethickness{0.15mm}
\bezier{100}(16.8, 2)(16.5, 4.2)(16.3, 4.5)
\bezier{100}(16.3, 4.5)(15.9, 5.2)(15.5, 5.2)
\bezier{100}(15.5, 5.2)(15, 5.2)(15, 4.5)
\bezier{100}(15, 4.5)(15, 3.8)(15.5, 3.8)
\bezier{100}(15.5, 3.8)(15.9, 3.8)(16.3, 4.5)
\bezier{100}(16.3, 4.5)(16.5, 4.8)(16.8, 7.5)

\bezier{100}(22.8, 2)(22.5, 4.2)(22.3, 4.5)
\bezier{5}(22.3, 4.5)(21.9, 5.2)(21.5, 5.2)
\bezier{5}(21.5, 5.2)(21, 5.2)(21, 4.5)
\bezier{5}(21, 4.5)(21, 3.8)(21.5, 3.8)
\bezier{5}(21.5, 3.8)(21.9, 3.8)(22.3, 4.5)
\bezier{100}(22.3, 4.5)(22.5, 4.8)(22.8, 7.5)
\end{picture}
\par
\noindent 
(2) $(x,y) \longmapsto (x,y^2,xy)$. The line $V(x)$ is melted with itself,
identifying the points $(0,1)$ and $(0,-1)$. $V(x) \longrightarrow V(r,t)$
is not birational, $C$ is a melting curve.
\par
\noindent
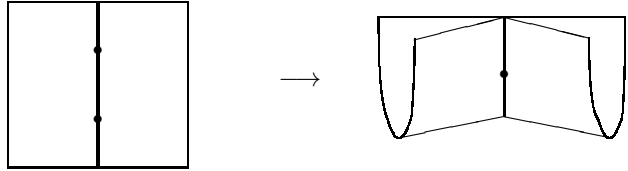
\unitlength0.4cm
\begin{picture}(30,8.5)

\put(5,3.6){\circle*{0.25}}
\put(5,5.9){\circle*{0.25}}

\put(2,2){\line(1,0){6}}
\put(2,2){\line(0,1){5.5}}
\put(2,7.5){\line(1,0){6}}
\put(8,2){\line(0,1){5.5}}
\thicklines
\put(5,2){\line(0,1){5.5}}

\put(11, 4.7){$\longrightarrow$}

\put(18.5, 5.1){\circle*{0.25}}
\thicklines
\put(18.5, 3.7){\line(0,1){3.3}}
\thinlines
\put(15, 3){\line(5,1){3.46}}
\put(18.5, 3.7){\line(5,-1){3.4}}

\put(14.3, 7){\line(1,0){4.2}}
\put(18.5, 7){\line(1,0){4}}
\put(15.5, 6.25){\line(4,1){3}}
\put(18.5, 7){\line(4,-1){2.8}}

\linethickness{0.15mm}

\bezier{100}(14.3, 7)(14.3, 4.3)(14.6, 3.6)
\bezier{100}(14.6, 3.6)(15, 2.3)(15.4, 3.7)
\bezier{100}(15.4, 3.7)(15.5, 5)(15.5, 6.25)

\bezier{100}(21.3, 6.3)(21.3, 4)(21.6, 3.6)
\bezier{100}(21.6, 3.6)(22, 2.3)(22.4, 3.7)
\bezier{100}(22.4, 3.7)(22.5, 4.6)(22.5, 7)

\end{picture}
\par
\noindent
(3) $(x,y) \longmapsto (x,y^2,y((y-1)^2+x^2)((y+1)^2+x^2),xy)$.
This identifies only the two points. $V(x)$ is birational with its image $C$,
thus $C$ is not a melting curve.
\par
\noindent
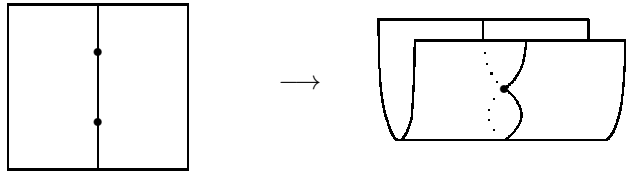
\unitlength0.4cm
\begin{picture}(30,8.5)

\put(5, 3.6){\circle*{0.25}}
\put(5, 5.9){\circle*{0.25}}

\put(2, 2){\line(1,0){6}}
\put(2, 2){\line(0,1){5.5}}
\put(2,7.5){\line(1,0){6}}
\put(8,2){\line(0,1){5.5}}
\put(5,2){\line(0,1){5.5}}

\put(11, 4.7){$\longrightarrow$}

\put(18.5,4.7){\circle*{0.25}}
\thicklines

\thinlines
\put(14.3, 7){\line(1,0){7}}
\put(15.5, 6.3){\line(1,0){7}}
\put(14.9, 3){\line(1,0){7}}

\linethickness{0.15mm}

\bezier{100}(14.3, 7)(14.3, 4.3)(14.6, 3.6)
\bezier{100}(14.6, 3.6)(15, 2.3)(15.4, 3.7)
\bezier{100}(15.4, 3.7)(15.5, 5)(15.5, 6.3)

\put(17.8,7){\line(0,-1){0.7}}
\bezier{50}(19.2, 6.3)(19.2, 5.2)(18.5, 4.7)
\bezier{5}(17.8, 6.3)(17.8, 5.2)(18.5, 4.7) 

\bezier{50}(18.5, 3)(19.6, 3.8)(18.5, 4.7)
\bezier{5}(18.5, 3)(17.4, 3.8)(18.5, 4.7)

\bezier{100}(21.9, 3)(22.5, 3.5)(22.5, 6.3)
\put(21.3, 7){\line(0,-1){0.7}}

\end{picture}
\par
\noindent
(4)
Consider the mapping $(x,y) \longmapsto (x,y^2,y(x^2-y^2(y+1))$
followed by the identification
of the points $(0,0,0)$ and $(-1,0,0)$.
Then $D=V(x^2-y^2(y+1)) \longmapsto C$ is not birational,
but $C$ (= the image of $D$) is not a melting curve for their common point.
Thus the necessary condition in prop. 3.4 (3) is not sufficient.

\par
\bigskip
\noindent
Holger Brenner
\par\noindent
Fakult\"at f\"ur Mathematik
\par\noindent
Ruhr-Universit\"at Bochum
\par\noindent
44780 Bochum
\par\noindent
brenner@cobra.ruhr-uni-bochum.de

\end{document}